\documentclass[12pt]{article}

\newtheorem{theorem}{Theorem}
\newtheorem{proposition}{Proposition}
\newtheorem{definition}{Definition}

\newtheorem{remark}{Remark}

\begin{document}

\author{Tornike Kadeishvili}


\title{Cochain operations defining Steenrod $\smile_i$-products in the bar construction}
\date{}

\maketitle

\begin{abstract}
The set of cochain multioperations defining Steenrod $\smile_i$-products in the
bar construction is constructed in terms of surjection operad. This structure
extends a Homotopy G-algebra structure which defines $\smile $ on the bar construction.
\end{abstract}

The Adams's cobar construction $\Omega C^*(X)$ of the chain complex of a
topological space $X$ determines the homology $H_*(\Omega X)$ of the loop
space just additively.

Lather on Baues \cite{Baues} had constructed the {\it geometric} diagonal
$$
\nabla _0:\Omega (C_{*}(X))\rightarrow \Omega (C_{*}(X))\otimes \Omega (C_{*}(X))
$$
which turns the cobar construction into a DG-Hopf algebra.
This diagonal allows to produce the next cobar construction
$\Omega \Omega (C_{*}(X))$ which models the double loop space.

Our aim here is to define on $\Omega (C_{*}(X))$ geometric
cooperations (dual to Steenrod's $\smile_i$-products)
$$
\{\nabla _i:\Omega (C_{*}(X))\rightarrow \Omega (C_{*}(X))\otimes
\Omega (C_{*}(X)),\ i=0,1,...\}
$$
satisfying the conditions
\begin{equation}
\label{steen}
deg\nabla _i=i,\ \nabla _id+(d\otimes 1+1\otimes d)\nabla _i=
\nabla _{i-1}+T\nabla_{i-1},
\end{equation}
(we work here on $Z_2$ so the signs are ignored).
These cooperations are necessary (but of course not sufficient)
for the further iteration of the cobar construction.

Namely we present particular elements $\{E^{i}_{p,q},\ i=0,1,...;\ p,q=1,2,...\}$
in the surjection operad (\cite{Mac}) $\chi $ such that the corresponding
chain multicooperations
$$
\{E^{i}_{p,q}:C_*(X)\to (C_*(X)^{\otimes p})\otimes
(C_*(X)^{\otimes q}),\ p,q=1,2,...\}
$$
define $\nabla _i $-s on the cobar construction $\Omega (C_{*}(X))$ and the corresponding
cochain multioperations
$$
\{E^{i}_{p,q}:(C^*(X)^{\otimes p})\otimes (C^*(X)^{\otimes q})
\to C^*(X),\ p,q=1,2,...\}
$$
define $\smile_i$-s on the bar construction $B(C^{*}(X))$.

It is known that the $\smile_i$-products in $C^*(X)$ are represented by the
following elements of $\chi $:
$$
\smile =(1,2);\ \smile_1 =(1,2,1);\ \smile_2 =(1,2,1,2);...\ .
$$
Let us consider them as first the line cochain operations.

The second line let be presented by a
{\it homotopy G-algebra} structure (\cite{Gerst-Vor}) on  $C^*(X)$ which consists of
the sequence of operations
$$
\{E_{1,q}:(C^*(X))\otimes (C^*(X)^{\otimes q})
\to C^*(X),\ q=1,2,...\},
$$
these operations in fact define a multiplication in the bar construction $BC^*(X)$. They
are represented by the following elements of $\chi $ (\cite{Mac}):
\begin{equation}
\label{1k}
E_{1,k}=(1,2,1,3,1,...,1,k,1,k+1,1).
\end{equation}

Bellow we present the next line cochain operations.
We introduce the notion of {\it extended homotopy G-algebra}, this is a
DG-algebra with certain additional structure which defines $\smile_i$-s on the bar
construction. A main example of such an object is again $C^*(X)$.
This structure consists of multioperations
$$
\{E^{i}_{p,q}:(C^*(X)^{\otimes p})\otimes (C^*(X)^{\otimes q})
\to C^*(X),\ i=0,1,...,\ p,q=1,2,...\}.
$$
We present particular elements $\{E^{i}_{p,q}\in \chi \}$ representing these
operations. Particularly, $E^0_{p,q}$ coincides with homotopy G-algebra structure (\ref{1k}),
$$
\begin{array}{ll}
E^1_{p,q}=&(1;p+1,1,p+2,1,...,p+q-1,1,p+q;\\
&1,p+q,2,p+q,3,...,p;p+q);
\end{array}
$$
and
$$
\begin{array}{ll}
E^2_{p,q}=\sum_{k=0}^{q-1}&(1;p+1,1,p+2,1,...,1,p+k+1;\\
&1,p+k+1,2,p+k+1,3,...,p+k+1,p;\\
&p+k+1,p,p+k+2,p,...,p+q;p).
\end{array}
$$

\section{ Homotopy G-algebras }

In this section we recall the notion of homotopy G-algebra from \cite{Gerst-Vor}
In order to extend it in the next section.

\subsection
{The notion of homotopy G-algebra}

\begin{definition}
A homotopy G algebra is a differential graded algebra (DG-algebra)
$(A,d,\cdot )$ together with a given sequence of multioperations
$$
E_{1,k}:A\otimes (\otimes ^kA)\rightarrow A,\quad
k=1,2,3,...\quad ,
$$
subject of the following conditions
\item
\label{HGAs}
 $deg E_{1,k}=-k,\  E_{1,0}=id ;$
\item
\begin{equation}
\label{HGAd}
\begin{array}{l}
d E_{1,k}(a;b_1,  ...,  b_k)+E_{1,k}(d a;b_1,  ...,
b_k)+
\sum_iE_{1k}(a;b_1,  ...,  d b_i,  ...,  b_k)= \\
b_1E_{1k}(a;b_2,  ...,  b_k)+\sum_iE_{1k}(a;b_1,  ...,
b_ib_{i+1},  ...,  b_k)+ \\
E_{1k}(a;b_1,  ...,  b_{k-1})b_k;
\end{array}
\end{equation}
\item
\begin{equation}
\label{HGAm}
\begin{array}{l}
a_1E_{1,k}(a_2;b_1,  ...,  b_k)+E_{1,k}(a_1\cdot a_2;b_1,
...,  b_k)+E_{1,k}(a_1;b_1,  ...,  b_k)a_2= \\
\sum_{p=1,...,k-1}E_{p,1}(a_1;b_1,  ...,  b_p)\cdot
E_{1,m-p}(a_2;b_{p+1},  ...,  b_k);
\end{array}
\end{equation}
\item
\begin{equation}
\label{HGAa}
\begin{array}{l}
E_{1,n} (E_{1,m}(a;b_1,...,b_m);c_1,...,c_n)=\\
\sum
E_{1,n-\sum l_i+m} (a;c_1,...,c_{k_1},
E_{1,l_1} (b_1;c_{k_1+1},...,c_{k_1+l_1}),
c_{ k_1+l_1+1},...,c_{k_m},\\
E_{1,l_m}(b_m;c_{k_m+1},...,c_{k_m+l_m}),
c_{ k_m+l_m +1},...,c_{n}).
\end{array}
\end{equation}
\end{definition}

Let us analyse these conditions in low dimensions.

For the operation $E_{1,1}$ the condition (\ref{HGAd}) gives
\begin{equation}
dE_{1,1}(a;b)+E_{1,1}(da;b)+ E_{1,1} (a;db)=a\cdot b + b\cdot a,
\end{equation}
i.e. the operation $E_{1,1}$ is sort of $\cup_1$ product, which measures the
noncommutativity of $A$.  Bellow we use the notation $E_{1,1}=\smile_1$.

The condition (\ref{HGAm}) gives
\begin{equation}
\label{Ehirsch}
(a\smile b)\smile_1 c+a\smile (b\smile_1c)+(a\smile_1c)\smile b=0,
\end{equation}
that is our $E_{1,1}=\smile_1$ satisfies so called {\it Hirsch formula},
which states that the map $f_b:A\to A$ defined as
$f_b(x)=x\smile_1b$ is a derivation.

The condition (\ref{HGAd}) gives
\begin{equation}
\label{Ehirschl}
\begin{array}{c}
a\smile_1( b\smile c)+b\smile (a\smile_1c)+(a\smile_1b)\smile c=\\
dE_{1,2}(a;b,c)+ E_{1,2}(da;b,c)+ E_{1,2}(a;db,c)+ E_{1,2}(a;b,dc),
\end{array}
\end{equation}
so the "left Hirsch formula" is satisfied just up to chain homotopy and
a homotopy is operation $E_{1,2}$, so this operation measures the lack of "left Hirsch formula".

Besides, the condition (\ref{HGAa}) gives
\begin{equation}
\label{Eassoc}
(a\smile_1 b)\smile_1c - a\smile_1 (b\smile_1 c)=E_{1,2}(a;b,c)+ E_{1,2}(a;c,b),
\end{equation}
so this $\smile_1$ is not strictly associative, but the operation $E_{1,2}$
somehow measures the lack of associativity too.


\subsection{Homotopy G-algebra structure and a multiplication in the bar construction}

For a homotopy G-algebra $(A,d,\cdot,\{E_{1,k}\})$
the sequence
$\{E_{1,k}\}$ defines in the bar construction $BA$ of a DG-algebra $(A,d,\cdot )$
a multiplication, turning $BA$ into a DG-Hopf algebra. In fact this means that a
homotopy G-algebra is a $B(\infty )$-algebra in the sense of \cite{Getzler}.

The sequence of operations $\{E_{1,k}\}$ defines a homomorphism
$$
E:BA\otimes BA\to A
$$
by
$E(a_1,...,a_m;b_1,...,b_n)=0$ if $m>1$ and $E(a;b_1,...,b_n)=E_{1,n}(a;b_1,...,b_n)$.

Since the bar construction $BA$ is a cofree coalgebra, a homomorphism $E$
induces a graded coalgebra map $\mu_E:BA\otimes BA \to BA$.

Then the conditions (\ref{HGAd}) and (\ref{HGAm}) are equivalent to the condition
$$
dE+E(d_{BA}\otimes id + id \otimes d_{BA})+E\smile E=0,
$$
that is $E$ is a twisting cochain, and this is equivalent to $\mu_E  $
being a chain map. Besides, the condition (\ref{HGAa}) is equivalent to
$\mu_E  $
being associative. Finally we have

\begin{proposition}
For a homotopy G-algebra $(A,d,\cdot,\{E_{1,k}\})$ the
bar construction $BA$ is a DG-Hopf algebra with respect to the standard
coproduct $\nabla_B:BA\to BA\otimes BA$ and the multiplication
$\mu _E:BA\otimes BA\to BA$.
\end{proposition}

\section{Extended homotopy G-algebras}

In this section we introduce the notion of {\it extended homotopy G-algebra}.
This is a DG-algebra with certain additional structure which defines
$\smile_i$-s on the bar construction.

\subsection{The notion of extended homotopy G-algebra}

\begin{definition}
An extended homotopy G-algebra we define as an object
$$
(A,d,\cdot,\{E^k_{p,q},\ k=0,1,...;\ p,q=1,2,...\})
$$
such that:
\item
$E^0_{p>1,q}=0$ and $(A,d,\cdot,\{E^0_{1,q}\})$ is a homotopy G-algebra;
\item
and
\begin{equation}
\label{ehga}
\begin{array}{c}
d E^{k}_{m,n}(a_1, ..., a_m;b_1, ..., b_n)+
\sum_iE^{k}_{m,n}(a_1, ..., da_i, ..., a_m;b_1,.., b_n)+\\
\sum_iE^{k}_{m,n}(a_1, ..., a_m;b_1, ..., db_i,.., b_n)+ \\
\sum_i
E^{k}_{m-1,n}(a_1, ..., a_i\cdot a_{i+1}, ...,a_m;b_1, ..., b_n)+\\
\sum_iE^{k}_{m,n-1}(a_1, ..., a_m;b_1, ..., b_i\cdot b_{i+1}, ..., b_n)+ \\
a_1E^{k}_{m-1,n}(a_2, ..., a_m;b_1, ..., b_n)+
E^{k}_{m-1,n}(a_1, ..., a_{m-1};b_1, ..., b_n)a_m \\
+b_1E^{k}_{m,n-1} (a_1, ..., a_m;b_2, ..., b_n)+
E^{k}_{m,n-1}(a_1, ..., a_m;b_1, ..., b_{n-1})b_n+ \\
\sum_{i=0 }^{k}\sum_{p,q}
T^iE^{k-i}_{p,q}(a_1, ..., a_p;b_1,.., b_q)\cdot
E^{i}_{m-p,n-q}(a_{p+1}, ..., a_m;b_{q+1}, ..., b_n)=\\
E^{k-1}_{m,n}(a_1, ..., a_m;b_1, ..., b_n)+
E^{k-1}_{n,m}(b_1, ..., b_n;a_1, ..., a_m),
\end{array}
\end{equation}
here $ TE^{i}_{p,q}(x_1, ..., x_p;y_1, ..., y_q)=
 E^{i}_{q,p}(y_1, ..., y_q;x_1, ..., x_p)$.
\end{definition}


Let us analyse this condition in low dimensions.

For the operation $E^k_{1,1}$ the condition (\ref{ehga}) gives
$$
dE^k_{1,1}(a;b)+E^k_{1,1} (da;b)+ E^k_{1,1} (a;db)= E^{k-1}_{1,1}(a;b)+ E^{k-1}_{1,1}(b;a),
$$
i.e. the operation $E^k_{1,1}$ is sort of $\smile_{k+1}$ product on $A$.
Bellow we use the notation $ E^k_{1,1}=\smile_{k+1}$.


Besides the condition (\ref{ehga}) also gives
\begin{equation}
\label{khirsch}
\begin{array}{c}
(a\smile b)\smile_k c+a\smile (b\smile_kc)+(a\smile_kc)\smile b+
E_{2,1}^{k-2}(a,b;c)+E_{1,2}^{k-2}(c;a,b)=\\
dE_{2,1}^{k-1}(a,b;c)+
E_{2,1}^{k-1}(da,b;c)+ E_{2,1}^{k-1}(a,db;c)+ E_{2,1}^{k-1}(a,b;dc)
\end{array}
\end{equation}
and
\begin{equation}
\label{khirschl}
\begin{array}{c}
a\smile_k( b\smile c)+b\smile (a\smile_kc)+(a\smile_kb)\smile c+
E_{1,2}^{k-2}(a;b,c)+ E_{2,1}^{k-2}(b,c;a)= \\
dE_{1,2}^{k-1}(a;b,c)+
E_{1,2}^{k-1}(a;db,c)+ E_{1,2}^{k-1}(a;b,dc),
\end{array}
\end{equation}
these are up to homotopy Hirsch type formulae connecting
$\smile_k$ and $\smile $. We remark here that a homotopy G-algebra
structure controls conection between $\smile $ and $\smile_1$, while
the extended homotopy G-algebra structure controls the connections
between $\smile $ and $\smile_k$-s (but not between $\smile_m $ and
$\smile_n$ generally).

As we already know a homotopy G-algebra structure defines a
multiplication in the bar construction. Bellow we are going to show that
an extended homotopy G-algebra structure defines on the bar construction
Steenrod $\smile_i$ products. But before we need some preliminary notions.

\subsection{DG-Hopf algebras with Seenrod coproducts}

Let a {\it $DG$-coalgebra with Steenrod's coproducts} be an object
$$
(A;d;\nabla _0,\nabla _1,\nabla _2,...\ )
$$
where $(A;d;\nabla _0)$ is a DG-coalgebra (with $\deg d=-1$),
and $\nabla _i:A\rightarrow A\otimes A,\ i>0,$ are cooperations, dual to
Steenrods $\smile _i$ products, i.e. they satisfy the conditions (\ref{steen})
$$
deg\nabla _i=i,\ \nabla _id+(d\otimes 1+1\otimes d)\nabla _i=
\nabla _{i-1}+T\nabla_{i-1}.
$$

Suppose now that $A$ additionally is equipped with a multiplication $\cdot :A\otimes A\to A$
which turns $(A,d,\cdot )$ into a $DG$-algebra. We are interested what kind
of compatibility of $\nabla _i$-s with the multiplication $\cdot $ should be
required.

The following notion was introduced in \cite{Kadender}, the dual notion was introduced by
V. Smirnov in \cite{Smir} and was called $\smile_{\infty }$-Hopf algebra.

\begin{definition}
A DG-Hopf algebra with Steenrod큦 coproducts we
define as an object
$$
(A,d,\cdot ,\nabla _0,\nabla _1,\nabla _2,...)
$$
where  $(A,d,\cdot )$
is a DG-algebra, $\nabla_i$-s satisfy (\ref{steen}) and additionally we
require the following connections between $\nabla _i$-s and
the product $\cdot $ ({\it decomposition rule}):
\begin{equation}
\label{decomp}
\nabla _n(a\cdot b)=\sum_{k=0}^n \nabla _{k}(a)\cdot T_{A\otimes A}^k \nabla _{n-k}(b),
\end{equation}
where $T_{A\otimes A}:A\otimes A\to A\otimes A$ is the permutation
map $T_{A\otimes A}(a\otimes b)= b\otimes a$ and $T^k$ it큦 iteration.
\end{definition}

Particularly $\nabla_0$ is a multiplicative map,
that is $(A,d,\cdot ,\nabla_0)$ is a DG-Hopf algebra; $\nabla_1$ is a
$(\nabla_0,T\nabla_0)$-derivation, etc.

This decomposition rule (\ref{decomp}) has the following sense: if $(A,\cdot )$
is a free (i.e tensor) algebra, that is $A=T(V)$ (for example the cobar construction),
then (\ref{decomp}) allows to construct the
cooperations $\nabla _i$ first on the generating vector space $V$ and then to
extend them on whole $A$ by the suitable {\it extension rule}, which
follows from the above decomposition rule (\ref{decomp}).

Let $(C,d,\Delta )$ be a DG-coalgebra and $\Omega C$ be its cobar
construction.
By definition $\Omega C$ is the tensor algebra $T(s^{-1}\bar{C})$
generated by the
desuspension $s^{-1}\bar{C}$ of the coaugmentation coideal $\bar{C}$.
So coproducts $\nabla_i: \Omega C\to \Omega C\otimes \Omega C$
satisfying (\ref{decomp}) are determined by their restrictions
$E^{i}:C\to \Omega C\otimes \Omega C$, which are homomorphis of degree $i-1$.

In order $\nabla_i$ to satisfy (\ref{steen}) $E^i$ should satisfy the condition
\begin{equation}
\label{twist}
\begin{array}{c}
d_{\Omega C\otimes \Omega C} E^i+E^i d+\sum_{k=0}^i
E^{k}\smile T^k_{\Omega C\otimes \Omega C}E^{i-k}=
E^{i-1}+T_{\Omega C\otimes \Omega C}E^{i-1},
\end{array}
\end{equation}
which is the restriction of (\ref{steen}) on $C$.

So if we want to construct on $\Omega C$ a sequence $\nabla_i$
forming a structure of DG-Hopf algebra with
Steenrod큦 coproducts we have to construct a sequence of
{\it higher twisting cochains} - homomorphisms $\{E^i,\ i=0,1,...;\ deg E^i=i-1\}$
satisfying (\ref{twist}). Note that $E^0$ is an ordinary twisting cochain:
$$
d_{\Omega C\otimes \Omega C} E^0+E^0 d+E^0\smile E^0=0.
$$

\subsection{DG-Hopf algebras with Seenrod products}

Here we dualise the previous section.

Let a {\it $DG$-algebra with Steenrod's products} be an object
$$
(A;d;\smile _0,\smile _1,\smile _2,...\ )
$$
where $(A;d;\smile _0)$ is a DG-algebra (with $\deg d=+1$),
and $\smile _i:A\otimes A\rightarrow A,\ i>0,$ are Steenrods $\smile _i$
products, i.e. they satisfy the conditions
\begin{equation}
\label{steenalg}
d(a\smile_ib)=da\smile_ib+a\smile_idb+a\smile_{i-1}b+b\smile_{i-1}b.
\end{equation}

Suppose now that $A$ additionally is equipped with a diagonal $\nabla:A\to A\otimes A$
which turns $(A,d,\nabla )$ into a $DG$-coalgebra. We are interested what kind
of compatibility of $\smile _i$-s with the diagonal $\nabla $
should be required.

\begin{definition}
A DG-Hopf algebra with Steenrod큦 products we
define as an object
$(A,d,\nabla ,\smile _0,\smile _1,\smile _2,...)$ where  $(A,d,\nabla )$
is a DG-coalgebra, the products $\smile_i:A\otimes A\to A$ satisfy
(\ref{steenalg}) and additionally we
require the following connections between $\smile _i$-s and
the diagonal $\nabla $:
\begin{equation}
\label{decompalg}
\nabla \cdot \smile_n =\sum_{k=0}^n (\smile_{k} \otimes \smile_{n-k}\cdot
T^k_{A\otimes A})\nabla_{A\otimes A}.
\end{equation}
\end{definition}

Particularly $\smile_0$ is a coalgebra map, i.e.
$(A,d,\nabla ,\smile_0)$ is a $DG$-Hopf algebra.


Let $(C,d,\cdot )$ be a DG-algebra and $BC$ be its bar
construction.
By definition $BC$ is the tensor coalgebra $T^c(s^{-1}\bar{C})$
generated by the
desuspension $s^{-1}\bar{C}$ of the augmentation ideal $\bar{C}$.

Since of cofreeness of $T^c$ products $\smile_i: BC\otimes BC\to BC$
satisfying (\ref{decompalg}) are determined by their projections
$E^{i}:BC\otimes BC\to BC \to C$, which are homomorphis of degree $1-i$.

In order $\smile_i$ to satisfy (\ref{steenalg}) $E^i$ should satisfy the condition
\begin{equation}
\label{twistalg}
\begin{array}{c}
dE^i+E^i (d_{BC}\otimes id+id\otimes d_{BC})+\sum_{k=0}^i
E^{k} \smile E^{i-k} T^k_{BC\otimes BC}=\\
E^{i-1}+ E^{i-1} T_{BC\otimes BC},
\end{array}
\end{equation}
which is the projection of (\ref{steenalg}) on $C$.

So if we want to construct on $BC$ a sequence $\smile_i$-s
forming a structure of DG-Hopf algebra with
Steenrod products we have to construct a sequence of
{\it higher twisting cochains} - homomorphisms $\{E^i,\ i=0,1,...;\ deg E^i=1-i\}$
satisfying (\ref{twistalg}). Note that $E^0$ is an ordinary twisting cochain:
$$
E^0d_{BC\otimes BC} E^0+dE^0 +E^0\smile E^0=0.
$$


\subsection{A structure of extended homotopy G-algebra and Steenrod products
in the bar construction}

As we already know the part of extended homotopy G-algebra - the
sequence of operations $\{E^0_{p,q}\}$ (which in fact is a homotopy
G-algebra structure) defines on the bar construction $BA$ a multiplication,
turning $BA$ into a DG-Hopf algebra.
Here we show that for an extended homotopy G-algebra
$(A,d,\cdot,\{E^k_{p,q}\})$
the sequence
$\{E^{k>0}_{p,q}\}$ defines in the bar construction $BA$ of a DG-algebra

$(A,d,\cdot )$ the $\smile_i$-products turning $BA$ into a
DG-Hopf algebra with $\smile_i$-s.

The sequences of operations $\{E^k_{p,q}\}$ define homomorphisms
$$
\{E^k:BA\otimes BA\to A,\ k=0,1,...\}
$$
by
$E^k(a_1,...,a_m;b_1,...,b_n)= E^k_{m,n}(a_1,...,a_m;b_1,...,b_n)$.

The condition (\ref{ehga}), which veirfy our $\{E^k_{p,q}\}$-s  is
equivalent to the condition (\ref{twistalg}) for the sequence $\{E^k\}$,
so they define the correct $\smile_k$-s on $BC$.

Finally we have
\begin{proposition}
For an extended homotopy G-algebra $(A,d,\cdot,\{E^k_{p,q}\})$ the
bar construction $BA$ is a DG-Hopf algebra with Steenrod $\smile_i$-products.
\end{proposition}


\section{Cochain complex $C^*(X)$ as an extended homotopy G-algebra}

Main example of an extended homotopy G-algebra gives the following
\begin{theorem}
The chain complex of a topological space $C_*(X)$ carries a structure of
extended homotopy G-algebra.
\end{theorem}

Bellow we construct particular
elements in the surjection operad $\chi $ (which acts on $C_*(X)$), which
represent the multicooperations
$\{ E^{k}_{p,q}\}$.


\subsection{Operations $E^k_{p,q}$ in the surjection operad}

Surjection operad $\chi $ \cite{Mac} is defined as a sequence
of chain complexes $\chi (n)$
where $\chi (n)_d$ is spanned by nondegenerate surjections
$u:(1,2,...,n+d)\to (1,2,...,n),\ u(i)\neq u(i+1)$.
For the structure maps of this operad, action of $\chi $ on $C_*(X)$ (on $C^*(X)$)
and the filtration
$ F_1\chi \subset...\subset F_n\chi \subset ... \subset \chi $, with $ F_n\chi $
equivalent to little n-cub operad, we refer to to \cite {Berger}.

A sirjection $u$ is written as a string $(u(1),u(2),..., u(n+d))$.
The string $(1,2)\in F_1\chi (2)_0$ corresponds to the cup-product $\smile $ in $C^*(X)$
(dually to the Alexander-Whitney diagonal $\Delta $ in $C_*(X)$);

the string $(1,2,1)\in F_2\chi (2)_1$ corresponds to $\smile_1$ in $C^*(X)$
(dually to $\Delta_1 $ in $C_*(X)$);
the string $(1,2,1,2)\in F_3\chi (2)_2$ corresponds to $\smile_2$ in $C^*(X)$
(dually to $\Delta_2 $ in $C_*(X)$) etc.

Here we present particular elements of $\chi $ representing operations
$E^k_{p,q}$. They are obtained from {\it admissible tables} which we define now.

The first row of an admissible table consists of single number $1$.
The second  looks as
$
p+1,1,p+2,1,...,1,p+k.
$
As we see it consists of the {\it stabile part} (1-s on even places)
and of the {\it increasing part} ($p+1,p+2,...,p+k$ on odd places).

Each next row starts with the
stabile number of the previous row which gives rise to increasing
part at odd places. As the stabile part at even places serves the
maximal number from the previous row. For example if the previous row ends with
$
..., i-2,s,i-1,s,i,
$
then the next row looks as
$
s,i,s+1,i,s+2,i,...,i,s+t.
$
Each row consists of odd number of terms.
Important remark: the stabile part of a row may be {\it empty},
in this case under the {\it stabile part} of this row we mean the maximal
number of the previous row.

As we see a table consists of two increasing sequences of integres
$1,2,...$ and $p+1,p+2,...$ (of course with repetitions and permutations).
The main restriction is that the first sequence should  necessarily end by $p$.
A table always ends with one term row.

After this we put all rows of the admissible table in one string and obtain an
{\it admissible string} in $\chi$. We say that this string belongs to $E^k_{p,q}$
if: it큦 table consists of  $k+3$ rows;  the first element of the second row
is $p+1$ and the maximal number which occurs in the string is $p+q$.

Here are some examples (by ; we indicate the ends of rows in admissible tables).
The admissible string
$$
(1;
5,1,6,1,7;
1,7,2,7,3,7,4;
7,4,8,4,9;
4)
$$
belongs to $E^2_{4,5}$ and
$$
(1;
4,1,5,1,6;
1,6,2;
6;
2,6,3;
6;
3)
$$
belongs to $E^4_{3,3}$.

An element $ E^k_{p,q}\in \chi $ we define as the sum of all
admissible strings belonging to it.

Particularly $E^{2k}_{1,1}=(1,2,1,...,1,2)$ and
$ E^{2k-1}_{1,1}=(1,2,1,...,1,2,1)$. They correspond to $\smile_{2k}$
and $\smile_{2k+1}$ respectively.

Besides
$ E^0_{1,q}=(1,2,1,3,...1,q+1,1)$.
These elements generate $F_2\chi $ (\cite {Mac}) and they determine on $C^*(X)$
a structure of homotopy G/-algebra..

Here are the examples of more higher operations:
$$
\begin{array}{ll}

E^1_{p,q}=&(1;p+1,1,p+2,1,...,p+q-1,1,p+q;\\
&1,p+q,2,p+q,3,...,p;p+q);
\end{array}
$$
$$
\begin{array}{ll}
E^2_{p,q}=\sum_{k=0}^{q-1}&(1;p+1,1,p+2,1,...,1,p+k+1;\\

&1,p+k+1,2,p+k+1,3,...,p+k+1,p;\\
&p+k+1,p,p+k+2,p,...,p+q;p).
\end{array}
$$

Generally $ E^k_{p,q}$ belong to filtration $F_{k+2}\chi $.


The way how these elements $ E^k_{p,q}$ are obtained is following.
In \cite{Kade-San2} it is shown that the bar construction $BC^*(X)$
actually is the cochain complex of certain cubical set, and in \cite{Kadender}
the Steenrod큦 $\smile_i$ products are constructed in the cochains of a cubical set.

\begin{remark}
The elements $ E^0_{1,q}$-s satisfy the defining conditions of a Homotopy G-algebra
already in $\chi $.
So do $ E^{k>0}_{p,q}$-s: they satisfy (\ref{ehga}) already in $\chi $. For example
the condition (\ref{Ehirsch}) is a result of
\begin{equation}
\label{hinchi}
(1,2,1)\circ_1(1,2)+(1,2)\circ_2(1,2,1)+(id\times T)(1,2)\circ_1(1,2,1)=0,
\end{equation}
(what is not a case in, say,
Barrat-Eccles operad: the suitable combination there is just
{\it homological to zero}. Note also that Barrat-Eccles operad
acts on $C(X)$ via $\chi $, see \cite{Berger}).
The condition (\ref{Ehirschl}) is a result of
$$
(1,2,1)\circ_2(1,2)+(T\times id)(1,2)\circ_2(1,2,1)+(1,2)\circ_1(1,2,1)=d(1,2,1,3,1).
$$
\end{remark}


\begin{remark}
The extended homotopy G-coalgebra structure, that
is the operations $E^k_{p,q}$, establishes connections just between
$\Delta_{k}$ and $\Delta$ (equivalently  between $\smile_k$ and $\smile $
or between $E^k_{1,1}$ and $1,2$ in $\chi $), but not a connections
between $\Delta_{m}$ and $\Delta_{n}$ generally. Here are two operations
establishing connections between $\Delta_{2}$ and $\Delta_{1}$:
$
G_{1,2}=(1,2,1,3,1,3,2),\  G_{2,1}=(1,2,3,2,3,1,3)\in F_3\chi (3)_4
$
satisfy the conditions
$$
\begin{array}{c}
dG_{2,1}(a,b;c)+ G_{2,1}(da,b;c)+ G_{2,1}(a,db;c)+ G_{2,1}(a,b;dc)=\\
(a\smile_1b)\smile_2c+ a\smile_1(b\smile_2c)+ (a\smile_2c)\smile_2b+
E^1_{2,1}(a,b;c)+E^1_{2,1}(b,a;c),
\end{array}
$$
and
$$
\begin{array}{c}
dG_{1,2}(a;b,c)+ G_{1,2}(da;b,c)+ G_{1,2}(a;db,c)+ G_{1,2}(a;b,dc)=\\
a\smile_2(b\smile_1c)+ b\smile_1(a\smile_2c)+ (a\smile_2c)\smile_1b+
E^1_{1,2}(a;b,c)+E^1_{1,2}(a;c,b),
\end{array}
$$
already in the operad $\chi $. Note that the element $(1,2)\in F_1\chi $ generates
the operad $F_1\chi $. furthermore, $(1,2)$ and
$E^1_{1,k}=(1,2,1,3,1,...,1,k,1,k+1,1)\in F_2\chi $
generate the operad $F_2\chi $ \cite{Mac} (but not freely: for example (\ref{hinchi})
is a relation). The elements  $G_{1,2}=(1,2,1,3,1,3,2),\  G_{2,1}=(1,2,3,2,3,1,3)\in
F_3\chi (3)_4$ should be a part of some rich structure, which, together with $(1,2)$
and $(1,2,1,3,...,1,k+1,1)$ generates $F_3\chi $.
\end{remark}

\begin{remark}
The elements $ E^0_{1,q}\in F_2\chi $ determine
on the Hochschild cochain complex of an associative algebra $C^*(A,A)$
operations from \cite{Kadehoch}, \cite{Getzler} and on the cochain complex $C^*(X)$
operations from \cite{Baues} forming in both cases a Homotopy G-algebra
structures \cite{Gerst-Vor}.
There is one more example where these
operations act, this is the cobar construction of a Hopf
algebra \cite{Kademeas}.
We remark also that the Hochschild cochain complex is not an extented homotopy G-algebra,
that is $ E^{k>0}_{p,q}\in F_3\chi $ do not act on $C^*(A,A)$ since of
nontriviality of Gerstehaber bracket in the Hochschild cohomology
(no $\smile_2$ in $C^*(A,A)$).
Here acts only the suboperad
$F_{2}\chi $ (Deligne conjecture), since it is generated by $(1,2) \in F_1\chi $ and
$ E^0_{1,q}\in F_2\chi $-s \cite{Mac}, whereas on $C^*(X)$ acts whole $\chi $.
\end{remark}

\vspace{1cm}
A. Razmadze Mathematical Institute of the Georgian Academy of Sciences,

M. Alxidze str. 1, Tbilisi, 380093,
Georgia

kade@rmi.acnet.ge

\end{document}